\documentclass[final,1p,times]{elsarticle}

\usepackage{amssymb,amsmath}
\usepackage{hyperref}
\usepackage[all]{xy}

\newtheorem{thm}{Theorem}[section]
\newtheorem{prop}[thm]{Proposition}
\newtheorem{lemma}[thm]{Lemma}
\newtheorem{cor}[thm]{Corollary}

\newproof{pf}{Proof}

\begin{document}


\begin{frontmatter}

\title{Four-class Skew-symmetric Association Schemes }

\author[JM]{Jianmin Ma}
\ead{Jianmin.Ma@emory.edu}
\author[KW]{Kaishun Wang\corref{cor}}
\ead{wangks@bnu.edu.cn}
\cortext[cor]{Corresponding author}
\address[JM]{Oxford College of Emory University, Oxford, GA 30054, USA}
\address[KW]{Sch. Math. Sci. \& Lab. Math. Com. Sys., Beijing Normal University, Beijing  100875, China}

\begin{abstract}
An association scheme is called skew-symmetric if it has no symmetric adjacency relations other than the diagonal one. In this paper, we investigate 4-class skew-symmetric association schemes. 
In recent work by the first author it was discovered that their character tables  fall into three types. We now determine their intersection matrices.  
We then determine the character tables for 4-class skew-symmetric pseudocyclic association schemes, the only known examples of which are cyclotomic schemes. As a result, we answer a question raised by S.Y. Song  in 1996.   We characterize and classify  4-class imprimitive skew-symmetric association schemes.  We also prove that none of 2-class Johnson schemes  admits a 4-class skew-symmetric fission scheme. Based on three types of character tables above, a  short list of feasible parameters is generated.   
\end{abstract}


\end{frontmatter}
\newcommand{\Cyc}{\mathrm{Cyc}}
\newcommand{\T}{\tt{T}}
\section{Introduction}
\noindent
 A {\em $d$-class association scheme} ${\mathfrak X}$ is a pair $(X, \{R_i\}_{i=0}^d)$, where $X$ is a finite set, and each $R_i$ is a nonempty subset of $X\times X$ satisfying following axioms:
\renewcommand{\theenumi}{\roman{enumi}}
\renewcommand{\labelenumi}{(\theenumi)}
\begin{enumerate}
\item \label{AS1}
$R_0 = \{(x,x)\;|\;x \in X\}$ is the diagonal relation;
\item \label{AS2} 
$X\times X=R_0\cup R_1\cup\ldots\cup R_d,\quad R_i\cap R_j=\emptyset$ $(i\ne j)$;
\item \label{AS3}
for each $i$, $R_i^{\T} = R_{i'}$ for some $0 \le  i' \le d$, where $R_i^{\T} = \{(y, x)\;|\; (x, y) \in  R_i\};$
\item \label{AS4}
there exists  integers  $p^k_{ij}$ such that for all $(x,y)\in R_k$ and all $i,j,k$,  
\[
p^k_{ij} = |\{z \in  X\;|\; (x, z) \in  R_i, (z, y)\in  R_j\}|.
\]
\end{enumerate}
 ${\mathfrak X}$ is also called an association  scheme with $d$ classes (a $d$-class association scheme, or even simply a scheme). The subsets $R_i$ are called the adjacency relations of ${\mathfrak X}$. The integers $p^k_{ij}$ are called the {\em intersection numbers} of ${\mathfrak X}$, and $k_i$ $(= p^0_{ii'})$ is called the \emph{valency} of $R_i$.
Furthermore, ${\mathfrak X}$ is called \emph{commutative} if 
 $p^k_{ij} =p^k_{ji} \mbox{ for all } i, j, k;$
${\mathfrak X}$ is called \emph{symmetric} if 
 ${R_i}^{\T} = R_i$  for all $i$.  
 
 In the rest of
this paper, all association schemes are assumed to be commutative. By a theorem of Higman \cite{Hig75}, association schemes with at most four classes are commutative. We refer readers to Bannai and Ito's book \cite{BaI84} for the general theory of association schemes. 

Let $ {\mathfrak X} = (X, \{R_i\}_{i=0}^d)$ be a commutative association scheme with $|X|=n$. The {\em adjacency matrix} $A_i$ of $R_i$ is the $n\times n$ matrix whose $(x,y)$-entry is 1 if  
$(x,y) \in R_i$ and 0 otherwise. By  the \emph{adjacency} or {\em Bose-Mesner algebra} $\mathfrak{A}$ of $\mathfrak{X}$ we mean the algebra generated 
by $A_0, \dots, A_d$ over the complex field. Axioms  
(\ref{AS1})-(\ref{AS4}) are equivalent to the following: 
\[
A_0 = I, \quad 
\sum_{i=0}^d A_i = J, \quad 
A_i^{\T} = A_{i'}, \quad
A_i A_j = \sum_{i=0}^d p^k_{ij} A_k,
\]
where $I$ and $J$ are the identity and all-one matrices of order $n$, respectively. 

Since $\mathfrak{A}$  consists of   commuting normal matrices, 
it has a second basis consisting of primitive idempotents $E_0= J/n, \dots, E_d$.  Let 
\[ E_i \circ E_j =\frac{1}{n} \sum_{k=0}^d q^k_{ij} E_k,
\] 
where $\circ$ is the Hadamard (entry-wise) product. The coefficients  $q^k_{ij}$ are called the 
\emph{Krein numbers} of $\mathfrak{X}$. The {\em Krein condition} asserts that all   $q^k_{ij}$ are 
nonnegative reals \cite[Theorem II.3.8]{BaI84}. 
The integers $m_i=\mathrm{rank}{E_i}$ are called the \emph{multiplicities} of $\mathfrak{X}$. In particular, $m_0=1$ is said to be  {\em trivial}. If $m_1=\dots= m_d$, $ {\mathfrak X}$ is called {\em pseudocyclic}. A pseudocyclic scheme  is equivalenced, i.e., all non-diagonal relations have the same valency (see \cite[p.48]{BCN} and references there).  

For $i=0,1, \dots, d$, let 
$$
  A_i = \sum_{j=0}^d  p_{ji} E_j.
$$
The following   $(d+1)\times (d+1)$ matrices is called the {\em character table}   of ${\mathfrak X}$:
$$
P = \begin{pmatrix} 
     1 & 1 & \dots & 1 \\
     1 & & \\
     \vdots & & {p_{ji}}& \\
     1 & & &
\end{pmatrix}, \quad {1\le j,i\le d.} 
$$

Both intersection numbers $p_{ij}^k$ and Krein numbers  $q_{ij}^k$ can be calculated from $P$ \cite[Section II.3]{BaI84}:
\begin{equation} \label {e:pijk}
p_{ij}^\ell = \frac{1}{n k_\ell} \sum_{h=0}^d m_h p_{hi} p_{hj}\overline{p}_{h\ell}, 
\quad 
q_{ij}^\ell = \frac{m_i m_j}{n } \sum_{h=0}^d\frac{1}{k_h^2} p_{ih} p_{jh}\overline{p}_{\ell h},
\end{equation}
where $\overline{p}_{h\ell}$ is the complex conjugate of 
${p}_{h\ell}$.  
 
For $i = 0,1,\dots, d$, the {\em $i$-th intersection matrix} $B_i$ is defined to be the $(d+1) \times (d+1)$ matrix whose $(j,k)$ entry is $p_{ij}^k$. The character table $P$ determines matrices $B_i$, and vice versa. 

One way to construct new association schemes is by merging  or splitting  relations in an existing scheme. More 
precisely, a partition $\Lambda_0, \Lambda_1,\ldots, \Lambda_e$ of the index set $\{0, 1,\ldots,d\}$ is
said to be \emph{admissible} \cite{ItM91} if $\Lambda_0=\{0\}, \Lambda_i \ne \emptyset$ and $\Lambda'_i = \Lambda_j$ for some 
$j\ (1\le i, j\le e)$, where $\Lambda_i^\prime = \{\alpha'\,| \alpha\in \Lambda_i\}$. Let
 $R_{\Lambda_i} = \cup_{\alpha \in \Lambda_i} R_{\alpha}$. If $\mathfrak{Y} =(X, \{R_{\Lambda_i}\}_{ i=0}^e)$ becomes an association 
scheme, it is called a \emph{fusion} scheme of ${\mathfrak X}$, while ${\mathfrak X}$ is called a \emph{fission} scheme of $\mathfrak{Y}$. If every admissible partition gives rise to a fusion 
scheme, ${\mathfrak X}$ is called {\em amorphous} \cite{ItM91}. 
  
In \cite{Ban93}, Bannai and Song raised a question regarding the existence of 4-class amorphous association schemes with {\em the diagonal relation being the only symmetric adjacency relation}. We showed the nonexistence of such schemes in  \cite{Ma08},  and this  implies the nonexistence of  amorphous association schemes with at least 4 classes with this property. These results bring our attention to what we call {\em skew-symmetric} association schemes, i.e.,   association schemes with no symmetric adjacency relations other than the diagonal one. In \cite{Ma08},  4-class skew-symmetric schemes are classified by their character tables, which  fall into three types.  
 In this paper, we investigate 4-class skew-symmetric association schemes.  
 
 The balance of this paper is structured as follows. We first determine the intersection matrices  of 4-class skew-symmetric association schemes (see Section \ref{s:ch}). We then determine the character tables  for 4-class skew-symmetric pseudocyclic association schemes. As a result,  we answer a question raised by Song \cite{Song96} about cyclotomic schemes (see Section \ref{s:pseu}).   In Section~\ref{s:imp}, we classify  4-class imprimitive skew-symmetric association schemes. In Section \ref{s:table}, we generate a short  list of feasible parameters. We conclude with some remarks in Section \ref{s:remark}.  

We conclude this section by briefly mentioning the status of association schemes with up to 4 classes.   Two-class skew-symmetric association schemes are regular tournaments.  Two-class symmetric association schemes (or equivalently,  strong regular graphs)  have been widely studied \cite{BCN, AEB2}. There are a few papers about 3-class symmetric 
association schemes \cite{Mat75,van99}. Three-class nonsymmetric schemes have been investigated  by Song \cite{Song95, Song96}, Goldbach and Claasen \cite{GoC96a, GoC96} and J{\o}rgensen \cite{Jorgensen}.  Though it is relatively easy to characterize the parameter sets of 3-class non-symmetric schemes, it is much hard to construct primitive examples \cite{Jorgensen}. Nonetheless, 
some families of 3-class non-symmetric schemes were constructed on Galois rings in characteristic 4 \cite {Lie88, Ma07}. In \cite{Song96}, Song initiated the study of 4-class  skew-symmetric association schemes.

\section{Character tables and Intersection matrices}\label{s:ch}
\noindent
Let $\mathfrak{X} = (X, \{R_0, R_1, R_2, R_2^{\T}, R_1^{\T}\})$ be a skew-symmetric association scheme.  So the symmetrization $\tilde{\mathfrak{X}} $ of  $\mathfrak{X}$ is a 2-class symmetric scheme, where  $\tilde{\mathfrak{X}} = (X, \{R_0, R_1 \cup R_1^{\T}, R_2\cup R_2^{\T}\})$. 

The concept of a strong regular graph is essentially the same as that of a 2-class symmetric association scheme. A regular graph $(X,R)$, with $n$ vertices and valency $k$, is called
 \emph{strongly regular} if any pair of adjacent vertices
have $\lambda$ common neighbors,  and any two distinct non-adjacent vertices have $\mu$
 common neighbors. We say this is an $(n, k, \lambda, \mu)$-strongly regular graph.
It is easy to verify that $\mathfrak{X} = (X, \{R_0, R, \bar{R}\})$ is a symmetric association scheme, 
where $(X,\bar{R})$ is the complement (graph) of $(X,R)$.  The parameters $n, k, \lambda, \mu$ of  $(X,R)$ determine those of  $\mathfrak{X}$.

The following result has been proved in \cite{Ma08}. 

\begin{thm} \label{t:rs20}
Let $\mathfrak{X} = (X, \{R_0, R_1, R_2, R_2^{\T}, R_1^{\T}\})$ be a skew-symmetric association scheme and 
let $\tilde{\mathfrak{X}} =( X, \{R_0, R_1\cup R_1^{\T} , R_2\cup R_2^{\T}\})$
be the symmetrization  of $\mathfrak{X}$. Let  
$$ \tilde{P} = \left[\begin{array}{ccc} 
      1 & k_1 & k_2 \\ 1  & r & t \\ 1 & s & u 
      \end{array}\right] 
\begin{array}{c} 1 \\ m_1 \\ m_2 \end{array} 
$$
be the character table of  $\tilde{\mathfrak{X}}$, where $1, m_1$ and $m_2$ are the multiplicities of the corresponding row entries (eigenvalues). The entries of $\tilde{P}$ and the multiplicities $m_i$
 can be calculated from the parameters 
$(n, k, \lambda, \mu)$ of the strongly regular graph $(X, R_1\cup R_1^{\T})$, where $k=k_1$. With a possible 
rearrangement of rows and columns, 
 the character table of $\mathfrak{X}$ has  the following form$:$
$$ 
P= \left[\begin{array}{ccccc} 
      1 & k_1/2  & k_2/2 & k_2/2 & k_1/2\\ 
      1  & \rho  & \tau & \bar{\tau} & \bar{\rho}\\ 
	 1 & \sigma  & \omega & \bar{\omega}& \bar{\sigma} \\
	 1 & \bar{\sigma} & \bar{\omega} & \omega & \sigma\\
     1 & \bar{\rho}  & \bar{\tau} & \tau & \rho 
      \end{array}\right]
 \begin{array} {c} 1 \\ m_1/2 \\ m_2/2 \\ m_2/2 \\ m_1/2 \end{array}.     
$$
 The entries $\rho,\omega, \tau,$ and $\sigma$  are one of the three cases:
\begin{enumerate}
\item \label{t:rs20i}
 $ \rho = r/2, \
   \sigma=  (s + \sqrt{-b})/2,\   
   \tau = (t + \sqrt{-z})/2, \
   \omega =  u/2,
$
where $  b = {nk_1}/{m_2},  z = {nk_2}/{m_1}.$
 \item  \label{t:rs20ii}
 $
 \rho =  \left( r + \sqrt{-y}\right)/2, \
  \sigma = {s}/{2}, \
\tau = {t}/{2}, \
\omega =\left( u + \sqrt{-c}\right)/2,
$
where 
$ y = {nk_1}/{m_1},  c ={nk_2}/{m_2}.
$
\item \label{t:rs20iii}
$ \rho = (r + \sqrt{-y})/2, \
   \tau =  (t + \sqrt{-z})/2, \
   \sigma=  (s + \sqrt{-b})/2, \
   \omega =  (u - \sqrt{-c})/2,
$
where all of $b, c, y,z$ are positive and  satisfy the following equations: 
\[
m_1 y + m_2b = nk_1, \quad m_1 z + m_2c = nk_2, \quad
 m_1\sqrt{yz} - m_2\sqrt{bc}=0 .
\]
If any of $y,z,b,c$ is taken to be a free variable,   the remaining  variables can be solved. Let $z$ be the free variable. We have
$$
b={\frac {m_1 zk_1}{k_2 m_2}},\quad
y={\frac {k_1  \left( nk_2-m_1 z \right) }{k_2 m_1}},\quad
c={\frac {nk_2-m_1 z}{m_2}}. 
$$
\end{enumerate}

\end{thm} 

We  denote by $P_I$, $P_{II}$ and $P_{III}$  the three character tables in the above theorem, and refer them as type I, II, and  III, respectively. 
For a skew-symmetric association scheme $\mathfrak{X} = (X, \{R_0, R_1, R_2, R_2^{\T}, R_1^{\T}\})$, the intersection matrices $B_1$ and $B_2$ determine $B_3$ and $B_4$.  For typographic convenience, we display only the principal part $B_i^{(0)}$ of an intersection matrix $B_i$, i.e., the lower-right 4 by 4 submatrix.
 
\begin{thm} \label{t:intn}
 With the notation in Theorem \ref{t:rs20}. Suppose that the strongly regular graph 
$(X, R_1\cup R_1^{\T})$ has parameters $n, k, \lambda, \mu$. Let $k=k_1$ and $k_2=n -k -1$. 

\begin{enumerate}
\item \label{t:intn1}
 For Theorem~\ref{t:rs20} {\rm (\ref{t:rs20i})}, we have
$$
B_1^{(0)} =
\left[\begin {array}{cccc} 
\frac {\lambda + s}{4}&
\frac {k(k-\lambda -1 -u)}{4k_2}&
\frac {k(k-\lambda -1 -u)}{4k_2}&
\frac {\lambda - 3s}{4}\\\noalign{\medskip}
\frac {k-\lambda -1 +u}{4} &
\frac {k-\mu + r}{4}&
\frac {k-\mu - r}{4}&
\frac {k-\lambda-1 - u}{4}\\\noalign{\medskip}
\frac {k-\lambda -1 +u}{4}&
\frac {k-\mu - r}{4}&
\frac {k-\mu + r}{4}&
\frac {k-\lambda-1 - u}{4}\\\noalign{\medskip}
\frac {\lambda + s}{4}&
\frac {k(k-\lambda -1 +u)}{4k_2}&
\frac {k(k-\lambda -1 +u)}{4k_2}&
\frac {\lambda + s}{4}
\end {array} \right], 
$$
$$
B_2^{(0)}=
\left[\begin {array}{cccc} 
\frac {k-\lambda -1 +u}{4}&
\frac {k-\mu + r}{4}&
\frac {k-\mu - r}{4}&
\frac {k-\lambda-1 - u}{4}\\\noalign{\medskip}
\frac {k_2(k-\mu - r)}{4k}&
\frac {(n-2k+\mu -2) + t}{4}&
\frac {(n-2k+\mu -2) -3t}{4}&
\frac {k_2(k-\mu - r)}{4k}\\\noalign{\medskip}
\frac {k_2(k-\mu + r)}{4k}&
\frac {(n-2k+\mu -2) + t}{4}&
\frac {(n-2k+\mu -2) + t}{4}&
\frac {k_2(k-\mu + r)}{4k}\\\noalign{\medskip}
\frac {k-\lambda-1 - u}{4}&
\frac {k-\mu + r}{4}&
\frac {k-\mu - r}{4}&
\frac {k-\lambda -1 +u}{4}
\end {array} \right]. 
$$
\item  \label{t:intn2}
For  Theorem~\ref{t:rs20} {\rm  (\ref{t:rs20ii})},
$$
B_1^{(0)} =\left[\begin {array}{cccc} 
\frac  {\lambda + r}{4}&
\frac  {k(k-\lambda -1 -t)}{4k_2}&
\frac  {k(k-\lambda -1 -t)}{4k_2}&
\frac  {\lambda - 3r}{4}\\\noalign{\medskip}
\frac  {k-\lambda -1 +t}{4} &
\frac  {k-\mu + s}{4}&
\frac  {k-\mu - s}{4}&
\frac  {k-\lambda-1 - t}{4}\\\noalign{\medskip}
\frac  {k-\lambda -1 +t}{4}&
\frac  {k-\mu - s}{4}&
\frac  {k-\mu + s}{4}&
\frac  {k-\lambda-1 - t}{4}\\\noalign{\medskip}
\frac  {\lambda + r}{4}&
\frac  {k(k-\lambda -1 +t)}{4k_2}&
\frac  {k(k-\lambda -1 +t)}{4k_2}&
\frac  {\lambda + r}{4}
\end {array} \right] ,
$$ 
$$
B_2^{(0)}=
\left[\begin {array}{cccc} 
\frac  {k-\lambda -1 +t}{4}&
\frac  {k-\mu + s}{4}&
\frac  {k-\mu - s}{4}&
\frac  {k-\lambda-1 - t}{4} \\\noalign{\medskip}
\frac  {k_2(k-\mu - s)}{4k}&
\frac  {(n-2k+\mu -2) +  u}{4}&
\frac  {(n-2k+\mu -2) -3u}{4}&
\frac  {k_2(k-\mu - s)}{4k}\\\noalign{\medskip}
\frac  {k_2(k-\mu + s)}{4k}&
\frac  {(n-2k+\mu -2) +  u}{4}&
\frac  {(n-2k+\mu -2) +  u}{4}&
\frac  {k_2(k-\mu + s)}{4k}\\\noalign{\medskip}
\frac  {k-\lambda-1 - t}{4}&
\frac  {k-\mu + s}{4}&
\frac  {k-\mu - s}{4}&
\frac  {k-\lambda -1 +t}{4}
\end {array} \right]. 
$$

\item \label{t:intn3}
For Theorem~\ref{t:rs20} {\rm (\ref{t:rs20iii})}, we have

$$
B_1^{(0)}=\left[\begin {array}{cccc} 
    {\frac   {nk\lambda +\Pi}{4nk}}&
    {\frac   { nk_2\mu + nk   + 2\Phi +\Pi  }{4nk_2}}&
    {\frac   { nk_2\mu + nk  -2\Phi +\Pi }{4nk_2}}&
    {\frac   {nk\lambda-3\Pi }{4nk}}\\\noalign{\medskip} 
    {\frac   {nk_2\mu -nk -\Pi  }{4nk}}&
    {\frac   { nk(n -2 k+\lambda) + \Gamma  }{4nk_2}}&
    {\frac   {nk(n -2 k+\lambda) -\Gamma + 2\Phi   }{4nk_2}}&
    { \frac   {  nk    +nk_2\mu -2\Phi+ \Pi }{4nk}} \\\noalign{\medskip} 
    {\frac   {nk_2 \mu -nk -\Pi }{4nk}}&
    {\frac   {nk(n -2 k+\lambda)-\Gamma -2\Phi  }{4nk_2}}&
    {\frac   {nk(n -2 k+\lambda) +\Gamma    }{4nk_2}}&
    {\frac   {  nk + nk_2 \mu +\Pi   +2\Phi  }{4nk}}\\\noalign{\medskip} 
    {\frac   {nk\lambda + \Pi}{4nk}  }&
    {\frac   {nk_2\mu-nk -\Pi  }{4nk_2}}&
    {\frac   {nk_2\mu-nk-\Pi  }{4nk_2}}&
    {\frac   { nk\lambda + \Pi}{4nk}}\end {array} \right] ,
$$
$$
B_2^{(0)}=\left[\begin {array}{cccc} 
    {\frac   {nk_2\mu  -nk -\Pi}{4nk}}&
    {\frac   {nk(n -2 k+\lambda) +\Gamma    }{4nk_2}}&
    {\frac   {nk(n -2 k+\lambda) -\Gamma + 2\Phi   }{4nk_2}}&
    { \frac   {  nk    +nk_2\mu -2\Phi+ \Pi }{4nk}}\\\noalign{\medskip} 
    {\frac   {nk(n -2 k+\lambda)-\Gamma -2\Phi }{4nk}}&
    {\frac   {nk_2(n -2 k+\mu) -\Gamma - 3nk_2 }{4nk_2}}&
    { \frac   {nk_2(n -2 k+\mu)+ nk_2 +3\Gamma }{4nk_2}}&
    {\frac   {nk(n -2 k+\lambda)+ 2\Phi -\Gamma  }{4nk}}\\\noalign{\medskip} 
    {\frac   {nk(n -2 k+\lambda) + \Gamma  }{4nk}}&
    {\frac   {nk_2(n -2 k+\mu) -\Gamma - 3nk_2 }{4nk_2}}&
    {\frac   {nk_2(n -2 k+\mu) -\Gamma - 3nk_2 }{4nk_2}}&
    {\frac   {nk(n -2 k+\lambda) + \Gamma }{4nk}}\\\noalign{\medskip} 
    {\frac   {  nk + nk_2 \mu +\Pi   +2\Phi  }{4nk}}&
    {\frac   {nk(n -2 k+\lambda) +\Gamma    }{4nk_2}}&
    {\frac   {nk(n -2 k+\lambda)-\Gamma -2\Phi  }{4nk_2}}&
    {\frac   {n k_2\mu  -nk -\Pi }{4nk}}
\end {array} \right] ,
$$
where 
$$
\Gamma = m_1 r z + m_2 s c, \quad \Phi = m_1 r \sqrt{yz} - m_2 s \sqrt{bc}, 
\quad \Pi = m_1 ry + m_2 bs.
$$  

\end{enumerate} 
\end{thm}
\begin{pf}
We now calculate the intersection number $p^1_{11}$ in (\ref{t:intn1}) as an example.
 
By Eq. (\ref{e:pijk}) and Theorem \ref{t:rs20} (\ref{t:rs20i}), we can obtain the (1,1)-entry 
$p_{11}^1$ of $B_1^{(0)}$:
\[
p_{11}^1 = {\dfrac {{{k_1}}^{3}+{ m_1}\,{r}^{3}+{s}^{3}{ m_2}+nsk_1}{4n{ k_1}}}.
\] 
Let $\tilde{A}_{0}, \tilde{A}_{1}, \tilde{A}_{2}$ and $\tilde{E}_0, \tilde{E}_1, \tilde{E}_2$ be the adjacency matrices and primitive idempotents of  $\tilde{\mathfrak{X}}$ in Theorem \ref{t:rs20}. Then
$\tilde{A}_1 = k_1 \tilde{E}_0 + r \tilde{E}_1 + s \tilde{E}_2$ and thus $\tilde{A}_1^3 = k_1^3 \tilde{E}_0 + r^3 \tilde{E}_1 + s^3 \tilde{E}_2$ under the 
usual matrix multiplication. Take the trace of both sides of this identity: 
\[\begin{array}{c}
\mathrm{trace} (k_1^3 \tilde{E}_0 + r^3 \tilde{E}_1 + s^3 \tilde{E}_2) = k_1^3 + m_1 r^3 + m_2 s^3, \\[3pt]
 \mathrm{trace}(\tilde{A}_1^3) = \mbox{sum of all entries of } (\tilde{A}_1^2 \circ \tilde{A}_1).
\end{array}\]
Since $\tilde{A}_1^2 = k_1 \tilde{A}_0 + \lambda \tilde{A}_1 + \mu \tilde{A}_2$,  $\tilde{A}_1^2 \circ \tilde{A}_1 = \lambda \tilde{A}_1$. Hence,
$\mathrm{trace}(\tilde{A}_1^3) = n k_1 \lambda$ and thus $ k_1^3 + m_1 r^3 + m_2 s^3 = n k_1 \lambda$.  
So $p_{11}^1 = (\lambda + s)/4$. The other intersection numbers can be 
obtained in a similar way.   \hfill $\Box$
\end{pf}

Since the character tables $P_I$ and $P_{II}$ are completely determined by $\tilde {P}$, we can check the 
parameter set $(n, k, \lambda ,\mu)$ of a strongly regular graph for possible fissions. Furthermore, 
that intersection numbers are nonnegative integers put some arithmetic restrictions on the parameters $n, k, \lambda ,\mu$.  Note the intersection matrices in Theorem \ref{t:intn} (\ref{t:intn2}) can be obtained from those in  (\ref{t:intn1}) by switching $r$ and $s$, and $t$ and $u$. Therefore, we will state these restrictions for the case (\ref{t:intn1}).

\begin{cor}\label{c:cor}
 For Theorem \ref{t:intn}  (\ref{t:intn1}), the following hold:
\begin{enumerate}
\item $r, s$ and hence $t, u$ must be integers;\label{c:cor1}
\item \label{c:cor2}
$\lambda + s \equiv 0 \bmod 4;$
\item \label{c:cor3}
$k (k -\lambda -1 + u) \equiv 0 \bmod 4k_2;$
\item \label{c:cor4}
$k_2 (k - \mu -r)  \equiv 0 \bmod 4k;$
\end{enumerate}
\end{cor}

\section{Pseudocyclicity }\label{s:pseu}

\noindent
A $d$-class association scheme is called \emph{pseudocyclic} if all the nontrivial multiplicities coincide.  Let $\mathrm{GF}(q)$ be a finite field with $q$ elements, and $\alpha$ be a primitive element of  $\mathrm{GF}(q)$. For a fixed divisor $d$ of $q-1$, 
define  
\[
(x,y)\in R_i \mbox{ if }  x-y \in \alpha^i \langle \alpha^d \rangle, \quad 1\le i \le d.
\]
These relations $R_i$ define a pseudocyclic association scheme. We denote it by $\Cyc(q,d)$, 
called the $d$-class  {\em cyclotomic scheme} over $\mathrm{GF}(q)$. 
The scheme $\Cyc(q,d)$ is symmetric if and only if  $(q-1)/d$ is even or $q$ is a power of 2. 

The intersection numbers of $\Cyc(q,d)$ are given by cyclotomic numbers of order $d$ 
(for example, see \cite[Chapter 2]{BEW98} for cyclotomy). The intersection number $p_{ij}^k$ is the number of 
nonzero elements $z\in  \mathrm{GF}(q)$ such $(x,z)\in R_i$ and $(z,y)\in R_j$ for any $(x,y) \in R_k$, which can be
reduced to  the number of solutions $s$ in $\alpha^{j-i}\langle \alpha^d \rangle$ such that 
$1 + s $ is in $\alpha^{k-i}\langle \alpha^d \rangle$. The latter is the cyclotomic number
$(j-i, k-i)$. 
It is not difficult to prove  the follow result. 
\begin{lemma}
Let $q=p^b$ a power of prime $p$. 
Then $Cyc(q,4)$ is a skew-symmetric scheme if and only if $q \equiv 5 \bmod 8$ and $b$ is  odd. 
\end{lemma}

The following result is due to Song \cite[Lemma 3.3]{Song96}.
\begin{lemma}
 Suppose that $P$ is the character table of a 4-class skew-symmetric scheme
whose symmetrization is $\Cyc(q,2)$ for $q\equiv 5 \bmod 8$.  Then $P$ has the following
form:
\[
P = \left[\begin{array}{ccccc}
1 & f & f & f & f  \\
1 & \rho & \bar {\rho} & \tau & \bar{\tau} \\
1 & \bar {\rho} & \rho & \bar{\tau} & \tau  \\
1 & \bar{\tau} & \tau  & \rho  & \bar {\rho} \\
1  & \tau & \bar{\tau} & \bar{\rho} & \rho  
\end{array}\right]
\begin{array}{c}  1 \\ f \\ f \\ f \\ f\end{array}, \quad
\mbox{ where $(*)$ }
\begin{cases}
 f = \frac{1}{4}(q-1)\\
 \rho + \bar{\rho} = \frac{1}{2} (-1 + \sqrt{q})\\
 \tau + \bar{\tau} = \frac{1}{2} (-1 - \sqrt{q})\\
 \rho\bar{\rho} + \tau\bar{\tau} = \frac{1}{8} (3q+1).
 \end{cases}
\]
\end{lemma}
Song \cite{Song96} gave two solutions for the system $(*)$ and he also raised the following questions:

\paragraph*{Question 1}  Is there any other solution for the system $(*)$ 
that yields a feasible fission table $P$ for given
$\tilde{P}$ (from the symmetrization of $P$), for a large prime power $q\equiv 5 \bmod {8}$?

\paragraph*{Question 2} 
Are there any other two-class primitive schemes (strongly regular graphs) that admit symmetrizable fission 
schemes besides $\Cyc(q,2)$ (Paley graph) for $q \equiv 5 \bmod 8$?\\

In the rest of this section, we will investigate these questions. 
For $q\equiv 1 \bmod 4$, the Paley graph $P(q)$ has as vertex set $\mathrm{GF}(q)$, with two vertices being adjacent if their difference is a nonzero square in $\mathrm{GF}(q)$. This graph gives rise to a strongly regular graph with $m_1= m_2 = (q-1)/2$.   
A strongly regular graph 
 with $m_1= m_2$ is  called a {\em conference graph}. Therefore, all Paley graphs are conference graphs. But the converse is not 
true; there are conference graphs that are not Paley graphs. A conference graph on $q$ vertices, denoted by 
$C(q)$, has parameters $k=(q-1)/2$, $ \lambda= (q-5)/4$ and $\mu= (q-1)/4$ ($q$ is not required to be a prime power). When there is no danger of confusion, we also use $P(q)$ and $C(q)$ for  the 2-class schemes from the Paley graph $P(q)$ and 
the conference graph $C(q)$, respectively.   For $C(q)$,  we determine the character tables  of  its 4-class skew-symmetric fission schemes.

\begin{thm} \label{p:ConfFission}
$C(q)$ has a putative skew-symmetric fission scheme with four classes if and only if $q\equiv 5 \bmod 8$,  
and there exist  integers $g, h$ such that $q = g^2 + 4h^2$ with $g\equiv 1 \bmod 4$. 
\end{thm} 
\begin{pf}

$C(q)$  has eigenvalues $ k=(q-1)/2$, $r, s = (-1 \pm \sqrt{q})/2$. 
We first show that $C(q)$  can not have  a skew-symmetric fission scheme with character table of type I or II. We treat  type II here and  type I  can be handled similarly. 
By Theorem \ref{t:intn} (\ref{t:intn2}), 
 \[
p^1_{11} = \frac{\lambda + r}{4} = \frac{q + 2\sqrt{q} -7 }{16}, \quad
p^2_{11} =   \frac  {k(k-\lambda -1 -t)}{4k_2} = \frac{q + 2 \sqrt{q} +1}{16}.
\]
Now $ p^2_{11} - p^1_{11} = 1/2$, a contradiction. 
 
Now suppose that $C(q)$ has a 4-class skew-symmetric fission scheme with character table of type III: 
\begin{equation} \label{e:table}
P_{III} = \left[ \begin {array}{lllll}
1 & f & f & f & f \\ 
1 &  \rho &    \tau & \bar{\tau} & \bar{\rho} \\
1 &  \tau & \bar{\rho} & \rho & \bar{\tau} \\
1 &  \bar{\tau} & \rho & \bar{\rho} & \tau \\
1 &  \bar{\rho} & \bar{\tau} & \tau & \rho 
\end {array} \right], 
\mbox{ where } 
\left\{ \begin{array} {l} 
f = \frac{1}{4}(q-1) \\ [5pt]
\rho = \frac{1}{4}\left(-1 + \sqrt{q} + 2 \sqrt{-(q-z)}\right) \\{\bigskip}
\tau = \frac{1}{4} \left(-1 - \sqrt{q} + 2 \sqrt{-z}\right).
\end{array}\right.
\end{equation}
We can calculate the first intersection matrices $B_1$ of the scheme from 
Theorem \ref{t:intn} (\ref{t:intn3}): 
\[
B_1 = \left[ \begin{array}{ccccc} 
    0 &  1 & 0 & 0 & 0 \\
    0 & A & B & D & C\\
    0 & E & E & B & D \\
    0 & E & D & E & B \\
    f & A & E & E & A
\end{array}\right],   \quad 
\mbox{ where } \quad
\left\{\begin{array}{l}
16q A =   q^2 - 7q + 2\sqrt{q}(q-2z)  \\
16q B = q^2 + q + 2\sqrt{q}(q-2z ) + 8\sqrt{q}\sqrt{z(q-z)}  \\
16q C =  q^2 +q -6\sqrt{q}(q-2z) \\
16q D =  q^2 + q + 2\sqrt{q}(q-2z) - 8\sqrt{q} \sqrt{z(q-z)} \\
16q E =   q^2 - 3q - 2\sqrt{q}(q-2z) .
\end{array}\right.
\]

Since $A+E = (q-5)/8$,  $q \equiv 5 \bmod 8$. So $q$ can not be a perfect square. Otherwise,  $\sqrt{q}$ is an odd integer and thus $\sqrt{q} \equiv 1$ or $3 \bmod 4$. It follows that $q \equiv 1 \bmod 8$, a contradiction.  

Now $q$ is a nonsquare integer. So $q - 2z = g\sqrt{q}$ for some integer $g$ in order for $A$  to be an 
integer. Since $A , C, E$ are integers, it follows that 
$g \equiv 1 \bmod 4$. In order for $B$ to be an integer, $\sqrt{z(q-z)} = h\sqrt{q}$ for some integer $h$.
Substituting $ z = (q-g\sqrt{q})/2$ into $\sqrt{z(q-z)} = h\sqrt{q}$, we can obtain that
\(
q = g^2 + 4 h^2.
\)

 Conversely, if $q\equiv 5 \bmod 8$  and if there exist  
integer $s, t$ such that $q = g^2 + 4h^2$ with $g\equiv 1 \bmod 4$, then  $P_{III}$ is the character table of  a putative scheme, where $ z = (q-g\sqrt{q})/2$. \hfill $\Box$ 
\end{pf}

We now determine the character table of the cyclotomic scheme $\Cyc(q,4)$ for $q \equiv 5\bmod 8$.

\begin{thm}\label{p:cycSch}
Let $q= 4f +1 $ be a prime power with $f$ odd. So $\Cyc(q,4)$ is a 4-class skew-symmetric association scheme. Let 
$g$ and $h$ be defined by $q= g^2 + 4h^2, g \equiv 1 \bmod {4}$ and $(g, q) = 1$; these conditions determine $g$ uniquely, and $h$ up to a sign.  Then the following hold: 
\begin{enumerate}
\item $\Cyc(q,4)$ has the following character table $P$: 
 \[
\left[ \begin {array}{lllll}
1 & f & f & f & f \\ 
1 &  \rho &    \tau & \bar{\tau} & \bar{\rho} \\
1 &  \tau & \bar{\rho} & \rho & \bar{\tau} \\
1 &  \bar{\tau} & \rho & \bar{\rho} & \tau \\
1 &  \bar{\rho} & \bar{\tau} & \tau & \rho 
\end {array} \right], 
\mbox{ where }
\begin{cases}
 f = \frac{1}{4}(q-1)\\
 \rho = \frac{1}{4} \left( -1 + \sqrt{q} + \sqrt{-2q - 2g \sqrt{q}} \right)
\\ {\bigskip}
\tau = \;\frac{1}{4} \left( -1 - \sqrt{q} + \sqrt{-2q + 2g \sqrt{q}} \right). 
 \end{cases}
\]
 
 \item \label{p:cycSch2}
 $\Cyc(q,4)$ has the following intersection matrices $B_1,B_2$:
\[
 \left[ \begin{array}{ccccc} 
    0 &  1 & 0 & 0 & 0 \\
    0 & A & B & D & C\\
    0 & E & E & B & D \\
    0 & E & D & E & B \\
    f & A & E & E & A
\end{array}\right],  \
 \left[ \begin {array}{ccccc} 
    0&0&1&0&0\\ 0 & E&E& B&D\\ 0& D &A &C & B \\ f &E& A & A & E
\\ 0&B&E&D&E\end {array} \right], 
\mbox {
where }
\left\{\begin{array}{l}
16A = q - 7 + 2g \\
16B = q  + 1 + 2g + 8h\\
16C = q + 1 - 6g\\
16D = q  + 1 +2g - 8h\\
16E = q - 3 - 2g.
\end{array}\right.
\]
\end{enumerate}
\end{thm}
\begin{pf}
Since the intersection numbers $p^k_{i\,j}$ of  $\Cyc(q,4)$ are given by the cyclotomic numbers $(j-i, k-i)$, which
are  given in  \cite[Proposition 11]{Myerson81}. From the intersection matrices, we can  
calculate the character table $P$. \hfill $\Box$
\end{pf}

We note that  a putative 4-class skew-symmetric fission scheme  of $C(q)$ has intersection matrices of the form in Theorem \ref{p:cycSch} (\ref{p:cycSch2}), where $q = g^2 + 4h^2$ with $g \equiv 1 \bmod 4$.
In particular,  Theorem \ref{p:ConfFission} answers Question 1, in which $q$ is a prime power with 
$q \equiv 5 \bmod 8$.  
Question 2 remains open. 
From  number theory (e.g., \cite[\S17.6]{Ire81}), $q = g^2 + 4h^2$ has a solution if and only if in the prime factorization of $q$,  every prime 
factor $\equiv 3 \bmod 4$ has an even exponent. In particular, if $q\equiv 5\bmod 8$ is a prime, there is essentially one way to express $q$ as a sum of two squares. So $\Cyc(q,2)$ has a unique fission table, which is realized by $\Cyc(q,4)$.


\section{Imprimitivity}\label{s:imp}

\noindent
In this section, we investigate 4-class skew-symmetric imprimitive association schemes. An association scheme $(X,\{R_i\}_{i=0}^d)$ is \emph{imprimitive} if the union of some relations is an equivalence relation that is not $R_0$ or $X\times X$, and {\em primitive} otherwise. 
 
Let $\mathfrak{X} = (X, \{R_0,R_1, R_2, R_2^{\T}, R_1^{\T}\}$ be a skew-symmetric association scheme. 
If $\mathfrak{X}$ is imprimitive, then $R_0 \cup R_1\cup R_1^{\T}$ or $R_0\cup R_2\cup R_2^{\T}$ is an equivalence relation on $X$. Without
loss of generality, we assume  that $R_0 \cup R_1\cup R_1^{\T}$ is an equivalence relation. So $(X, R_1\cup R_1^{\T})$ is the union of certain 
copies of a complete graph, $gK_f$, where $n=|X|= gf$. Thus, $(X, {R_2\cup R_2^{\T}})$ is the complete multipartite
graph $\overline { gK_f}$. The symmetrization $\tilde{\mathfrak{X}}$  has the following character table:
\begin{equation}\label{e:impCh}
\tilde{P} = 
\left[\begin{array}{ccc} 
 1 & f-1 & f(g-1)  \\ 1 & f -1 & -f \\ 1 & -1 & 0 
\end{array}\right]
\begin{array} {c} 1 \\ g-1 \\ g(f-1) \end{array}. 
\end{equation}
 
By Theorem \ref{t:rs20},  ${\mathfrak{X}}$ has one of the following character tables:
\[  
P_I = \left[\begin {array}{ccccc} 
1&\frac{f-1}{2} & \frac{f(g-1)}{2} & \frac{f(g-1)}{2} &\frac{(f-1)}{2} 
\\
1& \frac{f-1}{2}& \tau & \bar{\tau}& \frac{f-1}{2}
\\
1&\sigma &0&0& \bar{\sigma}
\\
1&\bar{\sigma} &0&0& {\sigma}
\\
1& \frac{f-1}{2}& \bar{\tau} & {\tau}& \frac{f-1}{2}
\end{array}\right],
\mbox{ where } 
\begin{cases}
\sigma = \dfrac{1}{2}\left(-1+ \sqrt{-f} \right)\\[5pt]
\tau =  \dfrac{1}{2}\left(-f + \sqrt {- gf^2}\right).
\end{cases}
\]
\[
P_{II} =  \left[\begin {array}{ccccc} 
1&\frac{f-1}{2} & \frac{f(g-1)}{2} & \frac{f(g-1)}{2} &\frac{(f-1)}{2} 
\\
1 & \rho & -\frac{f}{2} & -\frac{f}{2} & \bar{\rho} 
\\
1 & -1/2 & \omega & \bar{\omega} & -1/2 
\\
1 & -1/2 & \bar{\omega} & {\omega} & -1/2
\\
1 & \bar{\rho} & -\frac{f}{2} & -\frac{f}{2} & {\rho}
\end{array}\right], 
\mbox{where }
\begin{cases}
\rho = \dfrac{1}{2}\left(f-1+ \sqrt {- {\frac { fg ( f-1) }{g-1}}}\right) \\
\omega = \dfrac{1}{2}\sqrt {- {\frac {{f}^{2} ( g-1 ) }{f-1}}}.
\end{cases} 
\]
\[
P_{III} =  \left[\begin {array}{ccccc}
1&\frac{f-1}{2} & \frac{f(g-1)}{2} & \frac{f(g-1)}{2} &\frac{(f-1)}{2}  \\
1  & \rho  & \tau & \bar{\tau} & \bar{\rho}\\ 
	 1 & \sigma  & \omega & \bar{\omega}& \bar{\sigma} \\
	 1 & \bar{\sigma} & \bar{\omega} & \omega & \sigma\\
     1 & \bar{\rho}  & \bar{\tau} & \tau & \rho  
\end{array}\right],
\mbox{ where } 
\begin{cases}
\rho = \dfrac{1}{2}\left(f-1+\sqrt {-{\frac{ \left( {f}^{2}g-z \right) \left( f-1 \right) }{f \left( g-1 \right) }}}\right)\\[5pt]
\sigma = \dfrac{1}{2} \left( -1+\sqrt {-{\frac{z}{fg}}}\right)\\[5pt]
\tau = \dfrac{1}{2} \left(-f + \sqrt{-z} \right)\\[5pt]
\omega = \dfrac{-1}{2} \sqrt {-{\frac{ \left( {f}^{2}g-z \right)  \left( g-
1 \right) }{g \left( f-1 \right) }}}.
\end{cases}
\]

By Theorem \ref{t:intn} (\ref{t:intn1}), 
the intersection matrices $B_1, B_2$ for $P_I$ has the following principal parts:
\[ 
\dfrac{1}{4}
\left[ \begin {array}{cccc} 
f-3&0&0&f+1\\\noalign{\medskip}0&2(f -1) &0
&0\\\noalign{\medskip}0&0&2(f -1) &0\\\noalign{\medskip}f-3&0&0&f-3
\end {array} \right], 
 \dfrac{1}{4}\left[ \begin {array}{cccc} 
  0&2(f -1) &0&0 \\\noalign{\medskip}0&
 \left( g-3 \right) f& \left( g+1 \right) f&0\\\noalign{\medskip}2\,f
 \left( g-1 \right) & \left( g-3 \right) f& \left( g-3 \right) f&2\,f
 \left( g-1 \right) \\\noalign{\medskip}0&2(f-1) &0&0
\end {array}
 \right] . 
\]
It follows that  $f, g \equiv 3 \bmod 4$.   

The intersection matrix $B_1$ for $P_{II}$ has  negative entries and thus $P_{II}$ is not realizable.

For $P_{III}$,  the  the parameter $y$ can be calculate from Theorem \ref{t:rs20} (\ref{t:rs20iii}): $y = \frac{f-1}{f(g-1)} (gf^2 - z)$. So 
$z < gf^2$. By Theorem \ref{t:intn} (\ref{t:intn3}), we can obtain the intersection numbers. In particular, 
$p^1_{12} = \frac{z-gf^2 }{ 2fg}$ and hence $z \ge gf^2$, a  contradiction. So $P_{III}$ is not realizable either. 
Thus, we have proved following result, and  the ``\emph{if part}''  was mentioned in \cite[Lemma 2.3]{Song96}.

\begin{thm} \label{p:imp}
Let $\tilde{\mathfrak{X}}$ be an imprimitive 2-class association scheme with character table given in (\ref{e:impCh}).  
 $\tilde{\mathfrak{X}}$  has a putative 4-class skew-symmetric fission scheme if and only if  $f, g \equiv 3 \bmod 4$.
In this case, the character table of this fission scheme is realized as $P_I$.    
\end{thm}

As pointed out in \cite{Song96}, when $f, g \equiv 3 \bmod 4$ are prime powers, 
the fission scheme in Theorem \ref{p:imp} is realized as a wreath product of  $\Cyc(f,2)$ and $\Cyc(g,2)$.  For example, there are two 4-class skew-symmetric association schemes on 21 vertices from the wreath products 
$\Cyc(3,2)\; \mathrm {wr}\; \Cyc(7,2)$ and $\Cyc(7,2)\; \mathrm {wr}\; \Cyc(3,2)$.

\section{Lists of small feasible parameters}\label{s:table}

\noindent
In order to generate feasible parameters for  4-class skew-symmetric association schemes we shall classify them into
three sets according to their symmetrizations:
\begin{enumerate}
\item [1.] The symmetrization  is imprimitive;
\item [2.] 
The symmetrization is pseudocyclic thus a conference type, which is  primitive.  
\item [3.] 
The symmetrization is  non-conference type and primitive. 
\end{enumerate}

These three sets cover all possibilities. The tables to follow list information about existence and enumeration of primitive skew-symmetric schemes with 4 classes. In the heading \#,  a plus sign ``+''  means that there is at least  one such  scheme; a question mark ``?'' that there is no such scheme known;    and a number $m$ that there are $m$ such schemes.

Case 1 is determined by Theorem \ref{p:imp}, so each feasible parameter is determined by a pair $(f,g)$. For distinct primes $f$ and $g$ with $f, g\equiv 3 \bmod 4$, there are at least two 4-class skew-symmetric imprimitive association schemes on $fg$ vertices, realized by the  wreath product.  

 Case 2 is determined by Theorem \ref{p:ConfFission}. In this case,
$q = g^2 + 4h^2$ with $g \equiv 1 \bmod 4$, and the pairs $(g,h)$ can be easily generated (for example, the command \textbf{sum2sqr} in Maple).  Each  parameter set in this case  is determined by the pair $(g, h)$. In particular, if 
$q$ is a prime power with $(g, p) =1$, then $\Cyc(q, 4)$ realizes this parameter set. We list the parameter sets up to 325 vertices in Table \ref{table1}.  
\begin{table}[th]
\caption{ Feasible parameters of pseudocyclic skew-symmetric schemes}
\begin{center}
\begin {tabular}{c|*{13}r} 
$n$  & 5&13&29&37& 45&53&61&85&85&101&109&117&125\\\hline
$g$ & 1&-3&5&1&-3&-7&5&9&-7&1&-3&9&-11\\
\# & 1 &1 & 1&+&?&+&+&?&?&+&+&?&+\\\hline
$n$ & 125&149&157&173& 181&197&205&205&221&221&229&245&261\\
$g$ &5&-7 &-11&13& 9&1&-3&13&5&-11&-15&-7&-15\\
\# & ?&+ &+&+ & +&+&? &? &? & ? & + & ? & ?   \\\hline
$n$ & 269&277&293&317&325&325&325 & & & & & &  \\
$g$ & 13&9&17&-11&1&17&-15& & & & & & \\
\# & +& +  & +& + &?  &? &? & & & &  & &  
\end {tabular} 
\end{center}
\label{table1}
 \end{table}     

 An association 
scheme in Case 3 has character table given by one of the three types in Theorem \ref{t:rs20}. Character tables of type  I or  II are uniquely determined by the underlying  symmetrizations, while those of type III need an additional parameter $z$.  Here is how we determine $z$. Starting with the parameter set of a strongly regular graph, we generate the intersection matrices of all skew-symmetric fission schemes with four classes. Then we determine the value of $z$ using Theorem \ref{t:intn} (\ref{t:intn3}).   
Based on the table of strongly regular graphs on A. E. Brouwer's webpage \cite{AEB}, a list is generated for up to 1300 vertices in Table \ref{table2}. We further rule out some parameter sets with the Krein condition.

\begin{table}[th]
\caption{ Feasible parameters of non-conference type skew-symmetric schemes}
\begin{center}\begin{tabular}{llllllllll}
$n$ & $k$ & $\lambda$ & $\mu$ & $r^{m_1}$ & $s^{m_2}$ & type & $z$ & \# & \\\hline
57 & 14 & 1 & 4 & $2^{38}$ & $-5^{18}$ & III  & 27 & 0 & \cite{AEB} \\
105 & 26 & 13 & 4 & $11^{14}$ & $-2^{90}$ & III & 540 & 0 & Krein\\
253 & 42 & 21 & 4 & $19^{22}$ & $-2^{230}$ & III & 2300 & 0 & Krein \\
273 & 102 & 41 & 36& $11^{90}$ & $-6^{182} $& III & 364 & ? & \\
381 & 114 & 29 & 36 & $6^{254}$ & $-13^{126}$ & III & 147 & ?& \\
441 & 110 & 19  & 30 & $5^{530}$ & $-16^{110}$ & II  &  & ? & \\
441 & 110 & 19  & 30 & $5^{530}$ & $-16^{110}$ & III & 252 & ? &\\
465 & 58  & 29 & 4 & $27^{30}$ & $-2^{434}$ & III & 6076  & 0 & Krein\\
497 & 186 & 55 & 78 & $4^{426}$ & $-27^{70}$ & III & 175 & ? & \\
729 & 182 & 55  & 42 & $20^{182}$ & $-7^{546}$  & I & & ? & \\
741 & 74 & 37 & 4 & $35^{38}$ & $-2^{702}$ & III & 12636 & 0 & Krein \\
813 & 290 & 109 & 100 & $19^{270}$ & $-10^{542}$ & III & 1084 & ? & \\
889 & 222 & 35 &  62 & $5^{762}$ & $-32^{126}$ & II  & & ? & \\
889 & 222 & 35 &  62 & $5^{762}$ & $-32^{126}$ & III & 252& ? &  \\
945 & 354  & 153 & 120  & $39^{118} $& $-6^{826} $& II & & ? & \\
993 & 310 & 89 & 100 & $10^{662}$ & $-21^{330} $& III & 363& ? & \\
1065 & 266 & 103 & 54  & $53^{70}$ & $-4^{994}$ & III & 10224& ? &\\
1081 & 90 & 45 & 4& $43^{46}$ & $-2^{1034}$ & III & 22748&0 & Krein  \\
1225 & 306 & 89 & 72& $26^{306}$ & $-9^{918}$ & III & 1575& ? & \\
1225 & 510 & 215  & 210 & $20^{510}$ & $-15^{714}$ & I & & ? & \\
1241 & 310  & 81 & 76  & $18^{510}$ &  $-13^{730}$ & I& & ? & \\\hline
\end{tabular}\end{center}
\label{table2}
 \end{table}     

We note that the Krein condition rules out several parameter sets in Table \ref{table2}. It turns out that these parameter sets are from the family of  2-class Johnson schemes. The Johnson scheme $\mathrm{J}(v,2)$ has as vertices all 2-subsets of a $v$-set and two 2-subsets in 
relation $R_i$ if their intersection has size $2-i$ ($0\le i \le 2$). The relation $R_1$ defines a strongly regular graph, with parameters 
\[
n = \frac{v(v-1)}{2}, \quad k = 2(v-2), \quad \lambda = v-2, \quad \mu = 4.
\]
We ask the following:

\paragraph{Question 3} Does any Johnson scheme $\mathrm{J}(v,2)$ admit a 4-class skew-symmetric fission scheme?\\[-5pt]

We now settle  this question in the rest of this section. $\mathrm{J}(v,2)$ has the following character table:
\[
 P = \left[\begin{array} {ccc} 
      1 & 2(v-2) & \frac{1}{2} (v-2)(v-3) \\ 1 & v -4 & -v+3 \\ 1 & -2 & 1 
     \end{array}\right]
\begin{array} {c} 1 \\ v-1 \\ v(v-3)/2 \end{array}. 
\]
Since both of the nontrivial multiplicities $v-1$ and $v(v-3)/2$ are even, it is necessary that $v \equiv 3 \bmod 4$ for $\mathrm{J}(v,2)$ to admit a desired fission scheme. Since ${\lambda +s} = {v - 4} \equiv 3 \bmod {4}$, this fission scheme can not have 
character table of type I by Corollary \ref{c:cor} (\ref{c:cor2}). If its character table is  of type II, then by Theorem \ref{t:intn} (\ref{t:intn2}), we have $\frac{\lambda - 3r}{4} = \frac{5-v}{2}$ and 
$\frac{k-\mu + s}{4} = \frac{v-5}{2}$. So $v=5$, a contradiction.    
If this fission scheme has character table of type III,  we can compute its intersection numbers using Theorem \ref{t:intn} (\ref{t:intn3}). In particular, we have the following:
\[
\dfrac{nk\lambda + \Pi}{4nk } = \frac{v(v-3)^2 - 2z}{2v(v-3)},
 \quad 
\dfrac{nk(n -2 k+\lambda) +\Gamma}{4nk_2} = \frac{4z}{v(v-3)^2}.
\]
Since intersection numbers are nonnegative integer and $z > 0$, either $z = \frac{1}{4}v(v-3)^2$ or  $z = \frac{1}{2}v(v-3)^2$. The Krein numbers can be calculated from 
 Eq. (\ref{e:pijk}). For $z = \frac{1}{4}v(v-3)^2$, $q^3_{11} = \frac{(v-1)(v-1 - \sqrt{2v(v-1)})}{4(v-2)^2} < 0$,  a contradiction. 
For $z = \frac{1}{2}v(v-3)^2$, the value $c$ in Theorem \ref{t:rs20} (\ref{t:rs20iii}) is  
$ {\frac {nk_2-m_1 z}{m_2}} = \frac{1}{2}(v-1)(6-v)$. So $v = 3$ and hence $z = 0$, a contradiction. 
Therefore, we have proved:
\begin{prop}
No Johnson scheme $\mathrm{J}(v,2)$ admits a 4-class skew-symmetric fission scheme.
\end{prop}

\section {Concluding Remarks}\label{s:remark}

\paragraph{Remark 1}  We do not know   any association scheme with character table of type II or III.   
In view of Theorems \ref{p:ConfFission} and \ref{p:cycSch}, it is interesting to know which 4-class fission tables of $C(q)$ or $\Cyc(q,2)$  can be  realized. Using the small association schemes data of Hanaki and Miyamoto \cite{Hanaki}, we find that $\Cyc(q,4)$ are the only 4-class pseudocyclic skew-symmetric schemes  for $q\le 30$. 

\paragraph{Remark 2} 
The association schemes discussed in the paper are related to some combinatorial designs.   
Each pseudocyclic scheme give rise to a 2-$(n, f, f-1)$ design and 2-$(n, f+1, f+1)$ design(see \cite[p. 48]{BCN}). 
A 2-class skew-symmetric association scheme is equivalent to  a Hadamard tournament and an imprimitive skew-symmetric association scheme gives rise to a multipartite tournament (see \cite{Kirkland} ).  It seems interesting to pursuit further study along this line. 

\paragraph{Remark 3}  
Question 3 should be asked in more general form. Which 2-class association scheme admits a (putative) skew-symmetric fission scheme with 4 classes? For example, study the possibility of 4-class skew-symmetric fission scheme of known series of strongly regular graphs. We shall investigate this question in another paper.

\section*{Acknowledgments}
\noindent
The authors would like to thank Misha Muzychuk for helpful discussion and Qing Xiang for providing some references on cyclotomy. 
We are indebted to  the anonymous  reviewers  for their thoughtful reports,   and in particular, for  Remark 3 and  pointing out some mistakes  in an early draft.  K. Wang's research is supported by NCET-08-0052,  NSF of China (10871027) and the Fundamental Research  Funds for the Central Universities of China.

\end{document}